\documentclass[12pt]{article}
\usepackage{latexsym}
\usepackage{amsmath} 
\usepackage{amsfonts} 
\usepackage{amssymb} 
\usepackage{amsbsy} 
\usepackage{amsthm}
\usepackage[mathscr]{eucal} 

\def\qed{\relax\ifmmode\hskip2em \Box\else\nobreak\hskip1em $\Box$\fi}

\theoremstyle{plain}

\newenvironment{eq}{\begin{eqnarray*}}{\end{eqnarray*}}
\newenvironment{eqn}{\begin{eqnarray}}{\end{eqnarray}}

\newenvironment{notac}{\par\smallskip\noindent{\bfseries Notation}: }{\par\smallskip\normalfont}

\newenvironment{lst}{%
\begin{list}{--}{\setlength{\parsep}{1pt plus .5pt minus .2pt}%
                \setlength{\topsep}{3pt plus 1pt minus 2pt}
                \setlength{\itemsep}{3.0pt plus 1.5pt minus .6pt}}}%
{\end{list}}

\theoremstyle{definition}
\newtheorem{teo}{Theorem}[section]
\newtheorem{pro}[teo]{Proposition}
\newtheorem{lem}[teo]{Lemma}
\newtheorem{cor}[teo]{Corollary}
\newtheorem{den}[teo]{Definition}
\newtheorem{obs}[teo]{Remark}
\newtheorem{prp}[teo]{Propiedad}
\newtheorem{eje}[teo]{Example}
\newtheorem{prps}[teo]{Propiedades}
\newtheorem{obss}[teo]{Remarks}
\newtheorem{ejes}[teo]{Examples}
\newenvironment{teor}{\smallskip\begin{teo}}{\end{teo}\smallskip}
\newenvironment{prop}{\smallskip\begin{pro}}{\end{pro}\smallskip}
\newenvironment{lema}{\smallskip\begin{lem}}{\end{lem}\smallskip}

\newenvironment{defi}{\smallskip\begin{den}}{\end{den}\smallskip}
\newenvironment{obse}{\smallskip\begin{obs}}{\end{obs}\smallskip}

\newenvironment{ejem}{\smallskip\begin{eje}}{\end{eje}\smallskip}

\newenvironment{dem}{\par\smallskip\noindent {\textit{Proof\/}}:\\ }{$\hfill\square$\smallskip\par\smallskip}

\def\a{\alpha}
\def\b{\beta}
\def\ga{\gamma}

\def\balpha{\boldsymbol{\alpha}}
\def\bbeta{\boldsymbol{\beta}}
\def\bgamma{\boldsymbol{\gamma}}
\def\btita{\boldsymbol{\theta}}
\def\be{\boldsymbol{e}}
\def\bz{\boldsymbol{z}}
\def\bzeta{\boldsymbol{z}}
\def\bJ{\boldsymbol{J}}

\def\n{\mathbb{N}}      
\def\N{\mathbb{N}}
\def\R{\mathbb{R}}
\def\C{\mathbb{C}}      

\def\ds{\displaystyle}
\def\sst{\scriptstyle}

\begin{document}
\begin{center}
Quasi-analyticity in Carleman ultraholomorphic classes
\end{center}
\begin{center}
by
\end{center}
\begin{center}
ALBERTO LASTRA and JAVIER SANZ (Valladolid)
\end{center}
\par
{\bfseries{Abstract.}}
{\small We give a characterization for two different concepts of quasi-analyticity in Carleman ultraholomorphic classes of functions of several variables in polysectors. Also, working with strongly regular sequences, we establish generalizations of Watson's Lemma under an additional condition related to the growth index of the sequence.}

\section{Introduction}

\indent\indent In 1886 H. Poincar\'e put forward the concept of asymptotic expansion at~$0$ for holomorphic functions defined in an open sector in $\C$ with vertex at $0$. He intended to give an analytic meaning to formal power series solutions (in general, non convergent) of ordinary differential equations at irregular singular points.
With the modern formulation of his definition, the next statements turn out to be equivalent:
\begin{lst}
\item[(i)] A function $f$ admits asymptotic expansion in~$S$.
\item[(ii)] The derivatives of $f$ remain bounded in proper and bounded subsectors of $S$.
\end{lst}
Formal power series and asymptotic expansions of Gevrey type do constantly appear in the theory of algebraic
ordinary differential equations and of meromorphic, linear or not, systems of
ordinary differential equations at an irregular singular point.
Let $S$ be a sector with vertex at 0 in the Riemann surface of the logarithm $\mathcal{R}$, $\a\ge 1$, and $f$ a holomorphic function in $S$. The following are equivalent:
\begin{lst}
\item[(i)] $f$ admits Gevrey asymptotic expansion of order $\a\ge1$ in $S$ (we write $f\in\mathcal{W}_{\a}(S)$).
\item[(ii)] For every proper and bounded subsector $T$ of $S$, there exist constants $c,A>0$ such that $\sup_{z\in T}|f^{(p)}(z)|\le cA^{p}p!^{\a}$, $p\in\N_0:=\{0,1,2,\ldots\}$.
\end{lst}

It is easily deduced that the map sending $f\in\mathcal{W}_{\a}(S)$ to the sequence $(f^{(p)}(0))_{p\in\N_0}$, where $f^{(p)}(0):=\lim_{z\to0,z\in T}f^{(p)}(z)$, is well defined and linear. For all these facts we refer the reader to the book of W. Balser~\cite{balser}

Generalizing this situation, given a sequence of positive real numbers $\textbf{M}=(M_{p})_{p\in\N_{0}}$ and a sector
$S$ with vertex at $0$ in $\mathcal{R}$, we define
$\mathcal{A}_{\textbf{M}}(S)$ as the space of holomorphic functions $f$ defined in $S$ for which there exists $A>0$ (depending on $f$) such that $\sup_{p\in\N_{0},z\in S}\frac{|D^{p}f(z)|}{A^{p}p!M_{p}}<\infty,
$
and consider the Borel map $\mathcal{B}$ sending every function $f\in\mathcal{A}_{\textbf{M}}(S)$ to $(f^{(n)}(0))_{n\in\N_{0}}$. A classical problem is that of quasi-analyticity, that is, giving a characterization for the injectivity of $\mathcal{B}$. For Gevrey classes of order $\a>1$ (i.e. for
$M_n=n!^{\a-1}$, $n\in\N_0$), the result is classical and it is called Watson's Lemma~\cite{Watson}:
The class is quasi-analytic if, and only if, the opening of $S$ is greater than $\pi(\a-1)$. In 1966, B. I. Korenbljum~\cite{kor} solved the general problem, as we recall in Theorem~\ref{teo30}, in terms of the non
convergence of the logarithmic integral
\begin{equation}
\label{e354}
\int^{\infty}\frac{\log T_{\widetilde{\textbf{M}}}(r)}{r^{1+1/(\ga+1)}}\,dr,
\end{equation}
where $\widetilde{\textbf{M}}=(n!M_n)_{n\in\N_0}$, $T_{\widetilde{\textbf{M}}}$ is the Ostrowski's function associated to the sequence $\widetilde{\textbf{M}}$ (see~(\ref{defiTMr})), and the sector has opening $\ga\pi$.

In the case of several variables, H.~Majima in 1983~\cite{Majima1,Majima2} introduced the so-called strong asymptotic development. To every function $f$ admitting strong asymptotic development in a fixed polysector, he associates a unique family $\mathrm{TA}(f)$, named of strong asymptotic expansion of $f$, consisting of functions obtained, as in the one variable case, as limits of the derivatives of $f$ when some of its variables tend to zero (see~(\ref{defelefamder})). The elements of $\mathrm{TA}(f)$ admit strong asymptotic expansion in the corresponding polysector, and are linked by certain coherence conditions (see~(\ref{limcondcohe})).
This concept enjoys all the usual algebraic properties, and the equivalence we have firstly mentioned for the one variable case holds, as it was proved in~\cite{jesusisla} (a more accessible work is~\cite{HernandezSanz1}). Therefore, it is natural to consider, for a given sequence of positive real numbers $\textbf{M}=(M_{p})_{p\in\N_{0}}$ and a polysector $S$ in $\mathcal{R}^{n}$, the space $\mathcal{A}_{\textbf{M}}(S)$ of the holomorphic functions $f$ in $S$ such that there exists $A>0$ (depending on $f$) with
$$
\sup_{\balpha\in\N_{0}^{n},\bz\in S}\frac{|D^{\balpha}f(\bz)|}{A^{|\balpha|}|\balpha|!M_{|\balpha|}}<\infty.
$$
We write $\mathfrak{F}_{\textbf{M}}(S)$ for the set of all coherent families in $S$ (as explained in Section~\ref{seccionda}) and consider the maps
$$\mathcal{B}:\mathcal{A}_{\textbf{M}}(S)\longrightarrow \C^{\N_0^n}\quad\hbox{and}\quad\textrm{TA}:\mathcal{A}_{\textbf{M}}(S)\longrightarrow \mathfrak{F}_{\textbf{M}}(S),$$
where the first one is defined as $\mathcal{B}(f):=\big(D^{\balpha}f(\mathbf{0})\big)_{\balpha\in\N_0^n},
$ with
$$
D^{\balpha}f(\mathbf{0}):=\lim_{\bzeta\to {\bf 0},\bzeta\in S}
D^{\balpha}f(\bzeta)\in\textrm{TA}(f),\qquad \balpha\in\N_0^n.
$$
These maps are homomorphisms and we say $\mathcal{A}_{\textbf{M}}(S)$ is quasi-analytic (respectively, (s) quasi-analytic) whenever $\mathcal{B}$ (resp. $\textrm{TA}$) is injective.\par

In the several variables case, the works of P. Lelong~\cite{Lelong} and W. A. Groening~\cite{Groening} allowed J. A. Hern\'andez~\cite{jesusisla} to obtain (s) quasi-analyticity results for ultraholomorphic classes in polysectors whose elements $f$ admit estimates for their derivatives $D^{\balpha}f$ in terms of a multi-sequence $(M_{\balpha})_{\balpha\in\N_0^n}$. Along the same lines, the second author~\cite{Sanz1} proved two Watson's type lemmas concerning the (s) quasi-analyticity and the quasi-analyticity of the Gevrey classes considered by Y. Haraoka~\cite{Haraoka}. Our aim in this paper is to obtain the corresponding results for the general classes $\mathcal{A}_{\textbf{M}}(S)$ introduced above, and for both quasi-analyticity concepts.\par

After giving some notation, Section~\ref{seccionsfr} is devoted to strongly regular sequences, that will appear in some of the results in last section, and to the definition of the growth index related to them. We also present the function $T_{\textbf{M}}$. In Section~\ref{seccionda} we recall the theory of strong asymptotic expansions and introduce the ultraholomorphic classes of functions we will deal with. In Section~\ref{sectionpral}, Korenbljum's result will allow us to give a characterization of quasi-analyticity in several
variables, in one or the other sense and for arbitrary sequences $\textbf{M}$, in
terms of an integral similar to (\ref{e354}) in which the role of $\gamma$ is now played by
$\overline{\gamma}=\max\{\gamma_j:j=1,\ldots,n\}$ (Proposition~\ref{prop40}),
or $\underline{\gamma}=\min\{\gamma_j:j=1,\ldots,n\}$
(Proposition~\ref{propcaraccasianalit}), where $\ga_{j}\pi$ stands for the opening of $S_{j}$, $j\in\{1,..,n\}$, with $S=\prod_{k=1}^{n}S_{k}$.
Next, thanks to classical results by
S. Mandelbrojt~\cite{mandel}, we give a new sufficient condition for
(s) quasi-analyticity in Proposition~\ref{prop17}. When one works with strongly
regular sequences $\textbf{M}$, it is possible to establish generalizations of
Watson's Lemma (Propositions~\ref{prop18} and~\ref{proplemaWatsoncasianalit})
under the additional condition (\ref{e82}) related to the growth index.
We would like to point out that these results for strongly regular sequences are new even in dimension one, and they generalize previous results by J. Schmets and M. Valdivia~\cite{schmets} and by V. Thilliez~\cite{thilliez}.\par

\section{Notations}

\indent\indent $\N$ will stand for $\left\{1,2,\ldots\right\}$, and $\N_0=\N\cup\{0\}$. For $n\in\N$, we put $\mathcal{N}=\{1,2,\ldots,n\}$. If $J$ is a nonempty subset of $\mathcal{N}$,  $\#J$ denotes its cardinal number.

We will consider sectors in the Riemann surface of the logarithm $\mathcal{R}$ with vertex at $0$. Let $\theta>0$. We will write $S_{\theta}=\big\{z:|\arg z|<\frac{\theta\pi}{2}\big\}$, the sector of opening $\theta\pi$ and bisecting direction $d=0$.

Let $S$ be a sector. A proper subsector $T$ of $S$ is a sector such that $\overline{T}\setminus\left\{0\right\}\subseteq S$. If moreover $T$ is bounded, we say $T$ is a bounded proper subsector of $S$, and write $T\prec S$.

A polysector is a cartesian product $S=\prod_{j=1}^n S_j\subset\mathcal{R}^n$ of sectors. A polysector $T$ is a proper subpolysector of $S$ if $T=\prod_{j=1}^n T_j$ with $\overline{T}_j\setminus\left\{0\right\}\subseteq S_j$, $j=1,2,\ldots,n$. $T$ is bounded if each one of its factors is.\par
Given $\bzeta\in\mathcal{R}^n$, we write $\bzeta_J$ for the restriction of $\bzeta$ to $J$, regarding
$\bzeta$ as an element of $\mathcal{R}^{\mathcal{N}}$.\par
Let $J$ and $L$ be nonempty disjoint subsets of $\mathcal{N}$. For
$\bzeta_J\in\mathcal{R}^J$ and $\bzeta_L\in\mathcal{R}^L$, $(\bzeta_J,
\bzeta_L)$ represents the element of $\mathcal{R}^{J\cup L}$ satisfying
$
(\bzeta_J, \bzeta_L)_J=\bzeta_J
$,
$
(\bzeta_J, \bzeta_L)_L=\bzeta_L
$;
we also write
$J^{\prime}=\mathcal{N}\setminus J$, and for $j\in \mathcal{N}$ we use $j^{\prime}$ instead of $\{j\}^{\prime}$.
In particular, we shall use these conventions for
multi-indices.\par

For $\btita=(\theta_1,\ldots,\theta_n)\in(0,\infty)^n$, we write $S_{\btita}=\prod_{j=1}^nS_{\theta_j}$ and $S_{\btita_J}=\prod_{j\in J}S_{\theta_j}\subset\mathcal{R}^J$.\par

If $\boldsymbol{z}=(z_{1},z_{2},\ldots,z_{n})\in\mathcal{R}^{n}$, $\balpha=(\a_{1},\a_{2},\ldots,\a_{n})$, $\bbeta=(\b_{1},\b_{2},\ldots,\b_{n})\in\N_0^{n}$, we define:
$$\begin{array}{lll}
\boldsymbol{1}=(1,1,\ldots,1), & \be_j=(0,\ldots,\stackrel{j)}{1},\ldots,0), \\
|\balpha|=\a_{1}+\a_{2}+\ldots+\a_{n}, & \balpha!=\a_{1}!\a_{2}!\cdots\a_{n}!, \\
|\boldsymbol{z}^{\balpha}|=|\boldsymbol{z}|^{\balpha}=
|z_{1}|^{\a_{1}}|z_{2}|^{\a_{2}}\ldots|z_{n}|^{\a_{n}}, & \boldsymbol{z}^{\balpha}=z_{1}^{\a_{1}}z_{2}^{\a_{2}}\cdots z_{n}^{\a_{n}}, \\
\boldsymbol{D}^{\balpha}=\frac{\partial^{\balpha}}{\partial\boldsymbol{z}^{\balpha}}=
\frac{\partial^{|\balpha|}}{\partial z_{1}^{\a_{1}}\partial z_{2}^{\a_{2}}\ldots\partial z_{n}^{\a_{n}}}, & \balpha\le\bbeta\Leftrightarrow\a_{j}\le\b_{j},\ j\in\mathcal{N}.
\end{array}$$
For $\bJ\in\N_0^n$, we will frequently write $j=|\bJ|$.\par

\section{Preliminaries}
\subsection{Strongly regular sequences}\label{seccionsfr}

In what follows, $\textbf{M}=(M_p)_{p\in\N_0}$ will always stand for a sequence of
positive real numbers, and we will always assume that $M_0=1$. 
We say:\par
($\a_0$) $\textbf{M}$ is {\it logarithmically convex\/} if $M_{n}^{2}\le M_{n-1}M_{n+1}$ for
every $n\in\N$.\par
($\mu$) $\textbf{M}$ is {\it of moderate growth\/} if there exists $A>0$ such that
\begin{eqn}\label{moderategrowth}
M_{p+\ell}\le A^{p+\ell}M_{p}M_{\ell},\qquad p,\ell\in\N_0.
\end{eqn}%
\par($\gamma_1$) $\textbf{M}$ satisfies the {\it strong non-quasianalyticity condition\/} if there exists $B>0$ such that
\begin{eq}
\sum_{\ell\ge p}\frac{M_{\ell}}{(\ell+1)M_{\ell+1}}\le B\frac{M_{p}}{M_{p+1}},\qquad p\in\N_0.
\end{eq}%
$\textbf{M}$ is said to be {\it strongly regular\/} if it satisfies properties
$(\a_0)$, $(\mu)$ and $(\gamma_1)$.\par
Of course, for a strongly regular sequence $\textbf{M}$ the constants $A$ and $B$ above may be taken to be equal,
and they must be not less than~1. 

The measurable function $T_{\textbf{M}}:(0,\infty)\to(0,\infty]$, which firstly appeared in this context in a work by A. Ostrowski~\cite{Ostrowski}, is given by
\begin{equation}\label{defiTMr}
T_{\textbf{M}}(r)=\sup_{p\in\N_0}\frac{r^p}{M_p},\qquad r>0.
\end{equation}
Following V. Thilliez~\cite{thilliez}, we define next the growth index of a strongly regular sequence.
\begin{defi}\label{defi198}
Let $\textbf{M}=(M_{p})_{p\in\N_{0}}$ be a strongly regular sequence, $\ga>0$. We say $\textbf{M}$ satisfies property $\left(P_{\ga}\right)$  if there exist a sequence of real numbers $m'=(m'_{p})_{p\in\N}$ and a constant $a\ge1$ such that: (i) $a^{-1}M_{p}\le M_{p-1}m'_{p}\le aM_{p}$, $p\in\N$, and (ii) $\left((p+1)^{-\ga}m'_{p}\right)_{p\in\N}$ is increasing.
\end{defi}

\begin{prop}\label{lemthi1}(\cite{thilliez},Lemma 1.3.2). Let $\textbf{M}=(M_{p})_{p\in\N_{0}}$ be a strongly regular sequence. Then:
\begin{lst}
\item[(i)] There exists $\ga>0$ such that $(P_{\ga})$ is fulfilled and there also exists $a_{1}>0$ such that $a_{1}^{p}p!^{\ga}\le M_{p}$ for every $p\in\N_{0}$.
\item[(ii)] There exist $\delta>0$ and $a_{2}>0$ such that $M_{p}\le a_{2}^{p}p!^{\delta}$ for every $p\in\N_{0}$.
\end{lst}
\end{prop}

\begin{defi}
\label{defi212}
Let $\textbf{M}$ be a strongly regular sequence. The {\textit{growth index}} of $\textbf{M}$ is
$$\ga(\textbf{M})=\sup\{\ga\in\R:(P_{\ga})\hbox{ is fulfilled}\}.$$
\end{defi}
According to Proposition~\ref{lemthi1}, we have $\gamma(\textbf{M})\in(0,\infty)$.

\begin{ejem} For the Gevrey sequence of order $\a>0$, $\textbf{M}_{\a}=(p!^{\a})_{p\in\N_{0}}$, we have $\ga(\textbf{M}_{\a})=\a$.
\end{ejem}

\subsection{Strong asymptotic expansions and ultraholomorphic classes in polysectors}\label{seccionda}

Every definition and result in this section can be generalized by considering functions with values in a complex Banach space $B$, sequences of elements in $B$, and so on. However, for our purpose it will be sufficient to consider $B:=\C$, as it will be justified in Remark~\ref{obseclasecasianalitBanach}.

Let $n\in\N$, $n\ge 2$ and $S$ be a polysector in $\mathcal{R}^n$ with vertex at $\bf 0$. Taking into account the conventions adopted in the list of Notations, we give the following

\begin{defi}\label{defidaf}
We say a holomorphic function $f:S\to \C$ admits \textit{strong asymptotic development\/} in $S$ if there exists a family
$$
{\cal F}=\left\{\,f_{\balpha_J}:
\emptyset\neq J\subset \mathcal{N},\ \balpha_J\in\n_0^J\,\right\},
$$
where $f_{\balpha_J}$ is a holomorphic function defined in $S_{J'}$ whenever $J\neq \mathcal{N}$, and $f_{\balpha_J}\in \C$ if $J=\mathcal{N}$, in such a way that, if for every $\balpha\in\N_0^n$ we define the function
$$
\textrm{App}_{\balpha}({\cal F})(\bzeta):=
\!\sum_{\emptyset\neq J\subset \mathcal{N}}\!(-1)^{\#J+1}
\!\sum_{
\scriptstyle\bbeta_J\in\N_{0}^J\atop\scriptstyle
\bbeta_J\le\balpha_J-{\bf 1}_J}
\frac{f_{\bbeta_J}(\bzeta_{J'})}{\bbeta_J!}\bzeta_J^{\bbeta_J},\qquad \bzeta\in S,
$$
then for every proper and bounded subpolysector~$T$ of~$S$ and every
$\balpha\in\N_{0}^n$, there exists $c=c(\balpha,T)>0$ such that
\begin{equation}\label{cotadesaasinfuer}
\big|f(\bzeta)-\textrm{App}_{\balpha}({\cal F})(\bzeta)\big|\le c|\bzeta|^{\balpha},\qquad\bzeta\in T.
\end{equation}
\end{defi}
$\mathcal{F}$ is called the \textit{total family\/} of strong asymptotic development associated to $f$.
The map $\textrm{App}_{\balpha}({\cal F})$, which is holomorphic in $S$, is the \textit{approximant\/} of order $\balpha$
related to $\mathcal{F}$. We write $\mathcal{A}(S)$ for the space of holomorphic functions defined in $S$ that admit strong asymptotic development in $S$.\par

The next result is due to J. A. Hern\'andez~\cite{jesusisla} and it is based on a variant of Taylor's formula that appears in the work of Y. Haraoka~\cite{Haraoka}.
%
%
%
%
%
%
%
%
%

\begin{teor}\label{equivdafespas}
Let $f$ be a holomorphic function defined in $S$. The following statements are equivalent:
\begin{lst}
\item[(i)] $f$ admits strong asymptotic development in $S$.
\item[(ii)] For every $\balpha\in\N_0^n$ and $T\prec S$, we have
\begin{equation}\label{cotaderidaf}
Q_{\balpha,T}(f):=\sup_{\bzeta\in T}|D^{\balpha}f(\bzeta)|<\infty.
\end{equation}
\end{lst}
If properties (i) or (ii) are fulfilled, for every non empty subset $J$ of $\mathcal{N}$ and every $\balpha_J\in\N_{0}^J$ we have
\begin{eqn}\label{defelefamder}
f_{\balpha_J}(\bzeta_{J'})=
\lim_{\sst \bzeta_J\to {\bf 0}_J\atop\sst \bzeta_J\in T_J}
D^{(\balpha_J,{\bf 0}_{J'})}f(\bzeta),\quad \bzeta_{J'}\in S_{J'},
\end{eqn}%
for every $T_J\prec S_J$; the limit is uniform on every $T_{J'}\prec S_{J'}$ whenever $J\neq \mathcal{N}$, what implies that
$f_{\balpha_J}\in{\cal A}(S_{J'})$ (${\cal
A}(S_{\mathcal{N}'})$ is meant to be $\C$).
\end{teor}

\begin{obse}\label{obserelacotaderiapro}
According to~(\ref{defelefamder}), the family $\mathcal{F}$ in Definition~\ref{defidaf} is unique, and will be denoted by $\textrm{TA}(f)$, whereas the approximants will be $\textrm{App}_{\balpha}(f)$ from now on.\par

For future reference, we will make explicit the relationship between the estimates given in~(\ref{cotadesaasinfuer}) and the ones in~(\ref{cotaderidaf}). Given $T,T'$ polysectors, with $T\prec T'\prec S$, and $\balpha\in\N_0^n$, let us define
$$
P_{\balpha,T}(f):=\sup_{\bzeta\in T}\frac{| f(\bzeta)-\textrm{App}_{\balpha}(f)(\bzeta)|}{|\bzeta|^{\balpha}},
$$
then we have:
\begin{lst}
\item[(i)] $\ds P_{\balpha,T}(f)\le \frac{1}{\balpha!}Q_{\balpha,T}(f)$.
\item[(ii)] There exists a constant $A>0$, only depending on $T$ and $T'$, such that
$$
Q_{\balpha,T}(f)\le A^{|\balpha|}\balpha!P_{\balpha,T'}(f).
$$
\end{lst}
\end{obse}
According to the previous theorem, the space $\mathcal{A}(S)$ is stable under differentiation, what failed to hold for other concepts of asymptotic expansion in several variables (see \cite{GerardSibuya,HernandezSanz1}).
%
%
%
\begin{obse}
In case $n=1$, the concept of strong asymptotic development agrees with the usual one, and given $f\in\mathcal{A}(S)$, we have $\textrm{TA}(f)$ is reduced to the family of coefficients (except for the factorial numbers) of the formal power series of asymptotic expansion of $f$:
$$
\textrm{TA}(f)=\{\,a_m\in \C\colon \ m\in\n_0\,\},\qquad\textrm{with}\qquad f\sim\sum_{m=0}^{\infty}\frac{a_m}{m!}z^m.
$$
\end{obse}
%
%
%
%
%
%
\begin{prop}[Coherence conditions]\label{condcohe}
Let $f\in\mathcal{A}(S)$ and
$$
\textrm{TA}(f)=\left\{\,f_{\balpha_J}:\emptyset\neq J\subset \mathcal{N},\ \balpha_J\in\n_0^J\,\right\}
$$
be its associated total family. Then, for every pair of nonempty disjoint subsets $J$ and $L$ of $\mathcal{N}$,
every $\balpha_J\in\N_{0}^J$ and $\balpha_L\in\N_{0}^L$, and every $T_L\prec S_L$, we have
\begin{eqn}\label{limcondcohe}
\lim_{\sst \bzeta_L\to {\bf 0}\atop\sst \bzeta_L\in T_L}
D^{(\balpha_L,{\bf 0}_{(J\cup L)'})}f_{\balpha_J}(\bzeta_{J'})=
f_{(\balpha_J, \balpha_L)}(\bzeta_{(J\cup L)'});
\end{eqn}%
the limit is uniform in each~$T_{(J\cup L)'}\prec S_{(J\cup L)'}$ whenever $J\cup L\neq \mathcal{N}$.
\end{prop}%
>From the relations~(\ref{limcondcohe}) we immediately deduce that for every nonempty subset $J$ of $\mathcal{N}$ and every $\balpha_J\in\N_{0}^J$,
$$
\textrm{TA}(f_{\balpha_J})=\{\,f_{(\balpha_J,\bbeta_L)}:
\emptyset\neq L\subset J',\ \bbeta_L\in\N_{0}^L\,\}.
$$

\begin{defi}\label{defifamcohe}
We say a family
$$
{\cal F}=\{\,f_{\balpha_J}\in{\cal A}(S_{J'}):
\emptyset\neq J\subset \mathcal{N},\ \balpha_J\in\N_{0}^J\,\}
$$
is \textit{coherent\/} if it fulfills the conditions given in~(\ref{limcondcohe}).
\end{defi}

Given a polysector $S$, $\mathfrak{F}(S)$ will stand for the set of coherent families consisting of functions $f_{\balpha_J}\in{\cal A}(S_{J'})$, endowed with a vector structure in a natural way. We can consider the maps
$$
\mathcal{B}:\mathcal{A}(S)\longrightarrow \C^{\N_0^n}\quad\hbox{and}\quad\textrm{TA}:\mathcal{A}(S)\longrightarrow \mathfrak{F}(S),
$$
where the first one is given by
$$
\mathcal{B}(f):=\big(D^{\balpha}f(\mathbf{0})\big)_{\balpha\in\N_0^n},
$$
with
$$
D^{\balpha}f(\mathbf{0}):=\lim_{\sst \bzeta\to {\bf 0}\atop\sst \bzeta\in T\prec S}
D^{\balpha}f(\bzeta)=f_{\balpha}\in\textrm{TA}(f),\qquad \balpha\in\N_0^n.
$$
Without going into details, we mention that these maps are linear and are well behaved under differentiation.

\begin{obse}\label{notaprim}
Let $f\in\mathcal{A}(S)$. The \textit{first order family\/} associated to $f$ is given by
$$
\mathcal{B}_1(f):=\{\,f_{m_{\{j\}}}\in\mathcal{A}(S_{j'}): j\in \mathcal{N},\ m\in\n_0\,\}\subset\textrm{TA}(f).
$$
The first order family consists of the elements in the total family that depend on $n-1$ variables.
For the sake of simplicity, we will write $f_{jm}$ instead of $f_{m_{\{j\}}}$, $j\in \mathcal{N}$, $m\in\N_{0}$. As it can be seen in~\cite[Section\ 4]{GalindoSanz}, by virtue of the coherence conditions, the knowledge of $\mathcal{B}_1(f)$ is enough to determine $\textrm{TA}(f)$ uniquely. In fact, if $\mathcal{B}_{1}(f)$ consists of null functions, then the same holds for $\textrm{TA}(f)$.
\end{obse}
\begin{defi}
\label{deficlasesRoumieu}
Let $n\in\N$ and $S$ be a polysector in $\mathcal{R}^n$ with vertex at $\bf 0$. Given a sequence $\textbf{M}=(M_{p})_{p\in\N_{0}}$ and $A>0$, we define the class
$\mathcal{A}_{\textbf{M},A}(S)$ consisting of the holomorphic functions $f:S\to \C$ such that
\begin{equation}\label{e60}
\sup_{\bJ\in\N_{0}^n,\bz\in S}\frac{|D^{\bJ}f(\bz)|}{A^{j}j!M_{j}}<\infty.
\end{equation}
\end{defi}
According to Theorem~\ref{equivdafespas}, every element $f$ in $\mathcal{A}_{\textbf{M},A}(S)$ admits strong asymptotic development in $S$, in some sense ``uniform", because the limits or estimates in~(\ref{cotadesaasinfuer}), (\ref{defelefamder}) and (\ref{limcondcohe}) are valid in the whole corresponding polysector.\par
Obviously, the maps $\mathcal{B}$ and $\textrm{TA}$ can be restricted to these classes.

\section{Quasi-analyticity and generalizations of\\ Watson's Lemma}\label{sectionpral}
%
%
%
%
%
%
We aim at studying the injectivity of the maps $\mathcal{B}$ and $\mathcal{B}_1$, what justifies the introduction of two different concepts of quasi-analyticity.

\begin{defi}
Let $n\in\N$, $S$ be a (poly)sector in $\mathcal{R}^{n}$ and $\textbf{M}=(M_{p})_{p\in\N_{0}}$ be a sequence. We say that $\mathcal{A}_{\textbf{M}}(S)$ is \textit{(s) quasi-analytic} if the conditions:
\begin{lst}
\item[(i)] $f\in\mathcal{A}_{\textbf{M}}(S)$, and
\item[(ii)] every element in $\textrm{TA}(f)$ is null (or, equivalently, every function in the family $\mathcal{B}_1(f)$ is null),
\end{lst}
together imply that $f$ is null in $S$ (in other words, the class is (s) quasi-analytic if the map $\mathcal{B}_1$ restricted to the class is injective).\par

We say that $\mathcal{A}_{\textbf{M}}(S)$ is \textit{quasi-analytic} if the conditions:
\begin{lst}
\item[(i)] $f\in\mathcal{A}_{\textbf{M}}(S)$, and
\item[(ii)] $\mathcal{B}(f)$ is the null (multi-)sequence,
\end{lst}
together imply that $f$ is the null function in $S$.
\end{defi}

\begin{obse}\label{obseclasecasianalitBanach}
As we mentioned before, Section~\ref{seccionda} can be generalized in a natural way by considering functions with values in a general Banach space $B$ and sequences of elements in $B$. Let us denote by $\mathcal{A}_{\textbf{M}}(S,B)$ the generalization of the class $\mathcal{A}_{\textbf{M}}(S)$.
If $\mathcal{A}_{\textbf{M}}(S)$ is quasi-analytic (respectively, (s) quasi-analytic), so it is the class $\mathcal{A}_{\textbf{M}}(S,B)$ for any complex Banach space $B$, as it is easily deduced from the following facts:
\begin{lst}
\item[(i)] A map $f:S\to B$ is null if, and only if, we have $\varphi\circ f\equiv 0$ for every linear continuous functional $\varphi\in B'$.
\item[(ii)] If $f\in\mathcal{A}_{\textbf{M}}(S,B)$ and $\mathcal{B}(f)\equiv 0$ (resp., $\textrm{TA}(f)\equiv 0$), then for every linear continuous functional $\varphi\in B'$ we have $\varphi\circ f\in\mathcal{A}_{\textbf{M}}(S)$ and $\mathcal{B}(\varphi\circ f)\equiv 0$ (resp., $\textrm{TA}(\varphi\circ f)\equiv 0$).
\end{lst}
For this reason, the study of quasi-analyticity (respectively, (s) quasi-analyt\-icity) will be carried out only for ultraholomorphic classes of complex functions.
\end{obse}

\begin{obse}
It is worth saying that in the one variable case both definitions of quasi-analyticity coincide, because the families $\textrm{TA}(f)$ and $\mathcal{B}(f)$ are the same one. In the general case, quasi-analyticity implies (s) quasi-analyticity, since $\mathcal{B}(f)$ is a subfamily of $\textrm{TA}(f)$ for every $f\in\mathcal{A}_{\textbf{M}}(S)$.
\end{obse}

\begin{notac}
In this section, given a sequence $\textbf{M}=(M_{p})_{p\in\N_{0}}$, we will write $\widetilde{\textbf{M}}$ for the sequence $(p!M_{p})_{p\in\N_{0}}$. $T_{\textbf{M}}$ or $T_{\widetilde{\textbf{M}}}$ will be the maps defined according to~(\ref{defiTMr}). On the other hand, we will study the convergence of integrals
$$
\int^\infty g(r)\,dr,
$$
meaning that they are considered in intervals of the form $(r_0,\infty)$, the value $r_0>0$ being irrelevant. Finally, given an element $\bgamma=(\ga_1,\ldots,\ga_n)\in(0,\infty)^n$, we will write $S_{\bgamma}=\prod_{j=1}^{n}S_{\ga_{j}}=\{(z_{1},\ldots,z_{n})\in\mathcal{R}^{n}:|\arg(z_{j})|<\ga_{j}\pi,j\in\mathcal{N}\},$
$$
\underline{\ga}=\min\{\ga_j:j=1,\ldots,n\},\qquad \overline{\ga}=\max\{\ga_j:j=1,\ldots,n\}.
$$
\end{notac}

The next theorem, which characterizes quasi-analytic classes in the one variable case, is due to~B.~I.~Korenbljum~\cite[Theorem\ 3]{kor}, and it is based on the classical results of~S.~Mandelbrojt~\cite{mandel}. Condition (iii) in the following theorem has been used by V. Thilliez in~\cite{thilliez}.

\begin{teor}
\label{teo30}
Let $\textbf{M}=(M_{p})_{p\in\N_{0}}$ be a sequence and $\ga>0$. The following statements are equivalent:
\begin{lst}
\item[(i)] The class $\mathcal{A}_{\textbf{M}}(S_{\ga})$ is quasi-analytic.
\item[(ii)] $\displaystyle\int^{\infty}\frac{\log T_{\widetilde{\textbf{M}}}(r)}{r^{1+1/(\ga+1)}}\,dr$ does not converge.
\end{lst}
If $\textbf{M}$ verifies $(\a_0)$, the following is also equivalent to the previous statements:
\begin{lst}
\item[(iii)] $\displaystyle\sum_{p=0}^{\infty}\Big(\frac{M_{p}}{(p+1)M_{p+1}}\Big)^{1/(\ga+1)}$ does not converge.
\end{lst}
\end{teor}

Thanks to this result, we can establish easy conditions equivalent to (s) quasi-analyticity in the several variables case.

\begin{prop}
\label{prop40}
Let $n\in\N$, $\textbf{M}=(M_{p})_{p\in\N_{0}}$ be a sequence that fulfills $(\a_{0})$, and $\bgamma=(\ga_{1},\ldots,\ga_{n})\in(0,\infty)^{n}$. The following statements are equivalent:
\begin{lst}
\item[(i)] The class $\mathcal{A}_{\textbf{M}}(S_{\bgamma})$ is (s) quasi-analytic.
\item[(ii)] There exists $j\in\{1,\ldots,n\}$ such that $\mathcal{A}_{\textbf{M}}(S_{\ga_{j}})$ is quasi-analytic.
\item[(iii)] $\displaystyle\int^{\infty}\frac{\log T_{\widetilde{\textbf{M}}}(r)}{r^{1+1/(\overline{\gamma}+1)}}\,dr$ does not converge.
\item[(iv)] $\displaystyle\sum_{p=0}^{\infty}\Big(\frac{M_{p}}{(p+1)M_{p+1}}\Big)^{1/(\overline{\gamma}+1)}$ does not converge.
\end{lst}
\end{prop}

\begin{dem}
The equivalence between (iii) and (iv) is immediate according to Theorem~\ref{teo30}.
In order to demonstrate that (iii) implies (ii) it is enough to observe that $\overline{\gamma}=\gamma_{j}$ for some $j\in\{1,\ldots,n\}$, and apply Theorem~\ref{teo30}. In the other direction, and due to the same theorem, if for some $j$ the class $\mathcal{A}_{\textbf{M}}(S_{\ga_{j}})$ is quasi-analytic, then
$$\displaystyle\int^{\infty}\frac{\log T_{\widetilde{\textbf{M}}}(r)}{r^{1+1/(\gamma_j+1)}}\,dr=\infty,
$$
and the same can be deduced for the integral in (iii) on observing that the integrand, positive for $r>1$, becomes greater when we substitute the value $\gamma_j$ by a greater one, $\overline{\gamma}$.

Now, we are going to proof that (iii) implies (i). We consider a function $f\in\mathcal{A}_{\textbf{M}}(S_{\bgamma})$ such that $\textrm{TA}(f)$ is the null family, and we will conclude $f$ is null. Let $j\in\{1,\ldots,n\}$ be an index for which $\overline{\gamma}=\ga_{j}$, so that Theorem~\ref{teo30} guarantees that $\mathcal{A}_{\textbf{M}}(S_{\gamma_j})$ is quasi-analytic. For every $\bz_{j'}\in S_{\bgamma_{j'}}$ we define the function $\tilde{f}_{\bz_{j'}}$ in $S_{\gamma_j}$ by $\tilde{f}_{\bz_{j'}}(z_{j})=f(z_{j},\bz_{j'})$. We easily deduce that $\tilde{f}_{\bz_{j'}}\in\mathcal{A}_{\textbf{M}}(S_{\gamma_j})$. Moreover, for every $m\in\N_{0}$ and every $T_{j}\prec S_{\gamma_j}$ we have
$$
\lim_{z_{j}\to0,z_{j}\in T_{j}}\tilde{f}_{\bz_{j'}}^{(m)}(z_{j})=\lim_{z_{j}\to0,z_{j}\in T_{j}}D^{me_{j}}f(z_{j},\bz_{j'})=f_{mj}(\bz_{j'})=0.
$$
So, $\mathcal{B}(f_{\bz_{j'}})$ is the null family, and we conclude that $f_{\bz_{j'}}$ is identically null. Varying  $\bz_{j'}$ in $S_{\ga_{j'}}$ we conclude that $f$ is null.

Finally, if we suppose condition (iii) is not fulfilled we are going to find a function $F\in\mathcal{A}_{\textbf{M}}(S_{\bgamma})$, not identically 0 and such that $\mathcal{B}_{1}(F)\equiv0$, so that $\mathcal{A}_{\textbf{M}}(S_{\bgamma})$ is not (s) quasi-analytic. Since we have
$$\int^{\infty}\frac{\log T_{\widetilde{\textbf{M}}}(r)}{r^{1+1/(\overline{\gamma}+1)}}\,dr<\infty$$
and $\ga_j\le\overline{\ga}$, it is obvious that
$$\int^{\infty}\frac{\log T_{\widetilde{\textbf{M}}}(r)}{r^{1+1/(\ga_{j}+1)}}\,dr<\infty\quad\hbox{ for every }j\in\{1,\ldots,n\}.$$
Applying Theorem~\ref{teo30}, for every $j$ we can guarantee the existence of a non-zero function $f_{j}\in\mathcal{A}_{\textbf{M}}(S_{\gamma_{j}})$ such that $\mathcal{B}(f_j)$ is null. Let $F$ be the map defined in $S_{\bgamma}$ by
$$F(\bz)=\prod_{j=1}^{n}f_{j}(z_{j}),\qquad \bz=(z_{1},\ldots,z_{n}).$$
It is clear that $F$ is holomorphic in $S_{\bgamma}$ and it is not null. Moreover, $F\in\mathcal{A}_{\textbf{M}}(S_{\bgamma})$. Indeed, taking into account that $\textbf{M}$ fulfills property $(\a_{0})$, for every $\balpha=(\a_{1},\ldots,\a_{n})\in\N_{0}^{n}$ we have
\begin{align*}
|D^{\a}F(\bz)|&=|f_{1}^{(\a_{1})}(z_{1})\cdot\ldots\cdot f_{n}^{(\a_{n})}(z_{n})|\\
&\le C_{1}A_{1}^{\a_{1}}\a_{1}!M_{\a_{1}}\cdot\ldots\cdot C_{n}A_{n}^{\a_{n}}\a_{n}!M_{\a_{n}}\\
&\le CA^{|\balpha|}|\balpha|!M_{|\balpha|},\quad \bz\in S_{\bgamma},
\end{align*}
for certain positive constants $C,C_{1},\ldots,C_{n},A,A_{1},\ldots,A_{n}$.\\ Let us put $\mathcal{B}_1(F)=\{F_{jm}:j\in\mathcal{N},\ m\in\N_0\}$.
For $m\in\N_{0}$, $j\in\{1,\ldots,n\}$ and $\bz_{j'}\in S_{\bgamma_{j'}}$ we have
$$
F_{jm}(\bz_{j'})=\lim_{z_{j}\to0,z_{j}\in S_{\ga_j}}D^{me_{j}}F(\bz)=\prod_{\ell=1,\,\ell\neq j}^nf_{\ell}(z_{\ell})\lim_{z_{j}\to0,z_{j}\in S_{\ga_j}}f^{(m)}_{j}(z_{j})=0,$$
as desired.
\end{dem}

As it can be observed, the fact that $\textbf{M}$ fulfills $(\a_{0})$ is only used in the equivalence of (iii) and (iv), and in  (i)$\Rightarrow$(iii).

The following auxiliary result will allow us to give a new sufficient condition for (s) quasi-analyticity.

\begin{lema}
\label{lem76}
In the conditions of Proposition~\ref{prop40}, the following statements are equivalent:
\begin{lst}
\item[(i)] If $f$ is a function defined and holomorphic in $S_{\bgamma}$ such that
\begin{equation}
\label{e80} |f(\bz)|\le\frac{M_{|\balpha|}}{|\bz|^{\balpha}},\quad\hbox{ for every }\bz\in S_{\bgamma}\hbox{ and }\balpha\in \N_{0}^{n},
\end{equation}
then $f$ is null in $S_{\bgamma}$.
\item[(ii)] The integral
\begin{equation}
\label{e83}
\displaystyle\int^{\infty}\frac{\log T_{\textbf{M}}(r)}{r^{1+1/\overline{\gamma}}}\,dr
\end{equation}
does not converge.
\item[(iii)] The series $\displaystyle\sum_{p=0}^{\infty}\Big(\frac{M_{p}}{M_{p+1}}\Big)^{1/\overline{\gamma}}$  does not converge.
\end{lst}
\end{lema}

\begin{dem}
The equivalence between (ii) and (iii), whenever $\textbf{M}$ fulfills $(\a_0)$, can be found in~\cite[Theorem\ 2.4.III]{mandel}.

\noindent(i)$\Rightarrow$(ii) Let us suppose that the integral in (\ref{e83}) converges. For every index $j\in\{1,\ldots,n\}$ we then have that the integral obtained substituting $\overline{\gamma}$ by $\gamma_{j}$ in (\ref{e83}) is also convergent. Applying again~\cite[Theorem\ 2.4.III]{mandel}, for every $j$ we can guarantee the existence of a non-zero function $f_{j}$, holomorphic in $S_{1}=\{z:|\arg(z)|<\pi/2\}$ and such that
$$
|f_{j}(z)|\le M_{p}/|z|^{\ga_{j}p}\quad\textrm{ for every }p\in\N_{0},\ z\in S_{1}.
$$
We define the function $F$ by
$$
F(\bz)=f_{1}(z_{1}^{1/\ga_{1}})\cdot\ldots\cdot f_{n}(z_{n}^{1/\ga_{n}}), \qquad \bz=(z_{1},\ldots,z_{n})\in S_{\bgamma}.
$$
It is clear that $F$ is well defined, it is not null and it is holomorphic in $S_{\bgamma}$. Moreover, taking into account property $(\a_{0})$, for every $\balpha=(\a_{1},\ldots,\a_{n})\in\N_{0}^{n}$ and $\bz=(z_{1},\ldots,z_{n})\in S_{\bgamma}$ we have
$$
|F(\bz)|\le\frac{M_{\a_{1}}}{|z_{1}|^{\a_{1}}}\cdot\ldots\cdot\frac{M_{\a_{n}}}{|z_{n}|^{\a_{n}}}\le\frac{M_{|\balpha|}}{|\bz|^{\balpha}},
$$
concluding that (i) does not hold.\par\noindent
(ii)$\Rightarrow$(i) Let us consider a function $f$ holomorphic in $S_{\bgamma}$ and such that
$$|f(\bz)|\le M_{|\balpha|}/|\bz|^{\balpha}, \qquad\textrm{ for every }\bz\in S_{\bgamma}, \ \balpha\in\N_{0}^{n}.
$$
We can suppose, without loss of generality, that $\overline{\gamma}=\ga_{n}$.\\ For every $(z_{1},\ldots,z_{n-1})\in S_{\bgamma_{n'}}$, let $g$ be the map defined in $S_{1}$ by
$$
g(z_{n})=f(z_{1},\ldots,z_{n-1},z_{n}^{\ga_{n}}),
$$
which is a holomorphic map in $S_{1}$. For every $p\in\N_{0}$ we can apply (\ref{e80}) with $\balpha=(0,\ldots,0,p)$, obtaining that
$$|g(z_{n})|=|f(z_{1},\ldots,z_{n-1},z_{n}^{\ga_{n}})|\le\frac{M_{p}}{|z_{n}|^{\ga_{n}p}}.$$
Applying again~\cite[Theorem\ 2.4.III]{mandel}, we deduce that $g$ is the null function. As $(z_{1},\ldots,z_{n-1})\in S_{\ga_{n'}}$ was an arbitrary point, we conclude $f$ identically vanishes in $S_{\bgamma}$, as desired.
\end{dem}

\begin{prop}
\label{prop17}
Let $n\in\N$, $\textbf{M}=(M_{p})_{p\in\N_{0}}$ be a sequence, and  $\bgamma=(\ga_{1},\ldots,\ga_{n})\in(0,\infty)^{n}$. Then, we have:
\begin{lst}
\item[(i)] If $\displaystyle\int^{\infty}\frac{\log T_{\textbf{M}}(r)}{r^{1+1/\overline{\gamma}}}\,dr$ does not converge, $\mathcal{A}_{\textbf{M}}(S_{\bgamma})$ is (s) quasi-analytic.
\item[(ii)] If $\textbf{M}$ verifies $(\a_{0})$ and $\mathcal{A}_{\textbf{M}}(S_{\bgamma})$ is (s) quasi-analytic, then for every $\tilde{\gamma}>\overline{\gamma}$ we have $\displaystyle\sum_{p=0}^{\infty}\Big(\frac{M_{p}}{M_{p+1}}\Big)^{1/\tilde{\gamma}}$ and $\displaystyle\int^{\infty}\frac{\log T_{\textbf{M}}(r)}{r^{1+1/\tilde{\gamma}}}\,dr$ do not converge.\end{lst}
\end{prop}
\begin{dem} In both implications we will use the fact that, whenever $z\in S_{\gamma}$ for some $\gamma>0$, then also $1/z\in S_{\gamma}$.\par\noindent
(i) Let us consider $f\in\mathcal{A}_{\textbf{M}}(S_{\bgamma})$ with $\textrm{TA}(f)\equiv0$, so that all of its approximants are null. There exist constants $C,A>0$ such that
$$
|D^{\balpha}f(\bz)|\le CA^{|\balpha|}|\balpha|!M_{|\balpha|},\qquad\balpha\in\N_{0}^{n},\ \bz\in S_{\bgamma}.
$$
Taking into account Remark~\ref{obserelacotaderiapro}, we deduce that
\begin{equation}
\label{e129}
|f(\bz)|=|f(\bz)-\textrm{App}_{\balpha}(f)(\bz)|\le CA^{|\balpha|}\frac{|\balpha|!}{\balpha !}M_{|\balpha|}|\bz|^{\balpha},\qquad\bz\in S_{\bgamma},\ \balpha \in \N_{0}^{n}.
\end{equation}
If $\textbf{N}=(N_{p})_{p\in\N_{0}}$ is the sequence given by $N_{p}=C(nA)^{p}M_{p}$ for every $p\in\N_{0}$, then from (\ref{e129}) and from the inequality $|\balpha|!/\balpha!\le n^{|\balpha|}$, valid for all $\balpha\in\N_{0}^{n}$, we have that
$$
\Big|f\big(\frac{1}{\bz}\big)\Big|=\Big|f\big(\frac{1}{z_{1}},\ldots,\frac{1}{z_{n}}\big)\Big|\le\frac{N_{|\balpha|}}{|\bz|^{\balpha}}.
$$
It is easy to check that the integral in (\ref{e83}) and
$$
\displaystyle\int^{\infty}\frac{\log T_{\textbf{N}}(r)}{r^{1+1/\overline{\gamma}}}\,dr
$$
are simultaneously convergent or divergent. By our hypothesis, we deduce the second one does not converge. Applying Lemma~\ref{lem76} we get that the map $f(1/\bz)$ is null in $S_{\bgamma}$, as desired.\par\noindent
(ii) Now, suppose there exists $\widetilde{\gamma}>\overline{\gamma}$ such that the integral $\displaystyle\int^{\infty}\frac{\log T_{\textbf{M}}(r)}{r^{1+1/\widetilde{\gamma}}}\,dr$ is convergent. Lemma~\ref{lem76} tells us there is a holomorphic map $f$ in $S_{\tilde{\gamma}}$, not identically 0 and such that
$$
|f(z)|\le M_{p}/|z|^{p}\qquad\textrm{ for every }z\in S_{\widetilde{\gamma}},\ p\in\N_{0}.
$$
Now, if we define the map $F$ in $S_{(\widetilde{\gamma},\ldots,\widetilde{\gamma})}$ as
$$
F(\bz)=f(1/z_{1})\cdot\ldots\cdot f(1/z_{n}),\qquad\bz=(z_{1},\ldots,z_{n})\in S_{(\widetilde{\gamma},\ldots,\widetilde{\gamma})},
$$
then, due to property $(\a_0)$, we have
$$
|F(\bz)|\le M_{|\balpha|}|\bz|^{\balpha},\qquad\bz\in S_{(\widetilde{\gamma},\ldots,\widetilde{\gamma})},\ \balpha\in\N_{0}^{n},
$$
that is, $F$ admits strong asymptotic development in $S_{(\widetilde{\gamma},\ldots,\widetilde{\gamma})}$ and all the elements of $\textrm{TA}(F)$ are null. According to Remark~\ref{obserelacotaderiapro}, and taking into account that $S_{\bgamma}$ is a proper subpolysector of $S_{(\widetilde{\gamma},\ldots,\widetilde{\gamma})}$, we deduce that $G:=F|_{S_{\bgamma}}$ belongs to $\mathcal{A}_{\textbf{M}}(S_{\bgamma})$, it is not identically 0 and $\textrm{TA}(G)$ is the null family, so $\mathcal{A}_{\textbf{M}}(S_{\bgamma})$ is not (s) quasi-analytic.
\end{dem}

All the previous results are given in terms of a sequence $\textbf{M}$, subject, at most, to the condition $(\a_{0})$. The following ones deal with the case of a strongly regular sequence $\textbf{M}$. First of all, we will prove a result that, in the one variable setting, was already obtained by V.~Thilliez~\cite{thilliez}. Our proof is different from his and it also deals with the case of several variables.

\begin{prop}\label{propgammapequecasianal}
Let $n\in\N$, $\textbf{M}=(M_{p})_{p\in\N_{0}}$ be a strongly regular sequence, and $\bgamma=(\ga_1,\ldots,\ga_n)\in(0,\infty)^n$ such that $0<\overline{\ga}<\ga(\textbf{M})$.
Then, the class $\mathcal{A}_{\textbf{M}}(S_{\bgamma})$ is not (s) quasi-analytic.
\end{prop}

\begin{dem}
Let us take two real numbers $\widetilde{\ga}$ and $\eta$ such that $\overline{\ga}<\widetilde{\ga}<\eta<\ga(\textbf{M})$. We then have that property $(P_{\eta})$ is fulfilled, so, by Definition~\ref{defi198}, there exist a sequence $(m'_{p})_{p\in\N}$ and a constant $a\ge 1$ such that the sequence $\left((p+1)^{-\eta}m'_{p}\right)_{p\in\N}$ is increasing and $a^{-1}m_{p}\le m'_{p}\le am_{p}$ for every $p\in\N$, where $m_p=M_p/M_{p-1}$.
Therefore, the series
\begin{equation}\label{e145}\sum_{p=1}^{\infty}\Big(\frac{1}{m_{p}}\Big)^{1/\widetilde{\ga}}\end{equation}
is of the same nature, convergent or not, as
\begin{eqnarray}
\sum_{p=1}^{\infty}\Big(\frac{1}{m'_{p}}\Big)^{1/\widetilde{\ga}}&=&\sum_{p=1}^{\infty}\frac{1}{(p+1)^{\eta/\widetilde{\ga}}}
\frac{1}{(m'_{p}(p+1)^{-\eta})^{1/\widetilde{\ga}}}\\
&\le&\frac{1}{(m'_{1}2^{-\eta})^{1/\widetilde{\ga}}}\sum_{p=1}^{\infty}\frac{1}{(p+1)^{\eta/\widetilde{\ga}}}.
\end{eqnarray}
Since $\widetilde{\ga}<\eta$, this last series is convergent and also the series in (\ref{e145}) is. It is enough to apply the second item in the previous result to conclude.
\end{dem}

We now show that, under an additional hypothesis, the result in the previous proposition turns out to be an equivalence.
In this way, as it will be discussed in a later example, the classical Watson's Lemma, valid for Gevrey classes, is extended.

\begin{prop}[Generalization of Watson's Lemma]
\label{prop18}
Let $\textbf{M}$ be strongly regular and let us suppose that
\begin{equation}
\label{e82}
\sum_{n=0}^{\infty}\Big(\frac{M_{n}}{M_{n+1}}\Big)^{1/\ga(\textbf{M})}=\infty
\end{equation}
(or, in other words, $\ds\int^{\infty}\frac{\log
T_{\textbf{M}}(r)}{r^{1+1/\ga(\textbf{M})}}\,dr=\infty$). Let $n\in\N$ and $\bgamma\in(0,\infty)^{n}$. The following statements are equivalent:
\begin{lst}
\item[(i)] $\overline{\gamma}\ge\gamma(\textbf{M})$.
\item[(ii)] The class $\mathcal{A}_{\textbf{M}}(S_{\bgamma})$ is (s) quasi-analytic.
\end{lst}
\end{prop}

\begin{dem}
We only have to proof that (i) implies (ii). If~(\ref{e82}) is fulfilled and $\overline{\gamma}\ge\gamma(\textbf{M})$, then $\ds\int^{\infty}\frac{\log
T_{\textbf{M}}(r)}{r^{1+1/\overline{\ga}}}\,dr$ does not converge, and it is enough to apply (i) in Proposition~\ref{prop17}.
\end{dem}

\begin{ejem}\label{ejemclasese82}
For $\a>0$ and $\b\ge 0$ we consider the sequence $\textbf{M}=(M_{p})_{p\in\N_0}$ given by
$$
M_{p}=p!^{\a}\Big(\prod_{k=0}^{p}\log(e+k)\Big)^{\b},\qquad p\in\N_0.
$$
It is not hard to check that $\textbf{M}$
is strongly regular and $\ga(\textbf{M})=\a$. Moreover, $\textbf{M}$ fulfills condition (\ref{e82}) if, and only if, $\b\le\a$. It is important to observe that when $\beta=0$ we get the Gevrey sequences, $\textbf{M}_{\a}=(p!^{\a})_{p\in\N_0}$, and consequently for every $\a>0$ we have that $\textbf{M}_{\a}$ fulfills~(\ref{e82}). So, the previous result generalizes Watson's Lemma.
\end{ejem}

\begin{obse}\label{probabie}
It is an open problem to decide whether the condition $\overline{\gamma}\ge\gamma(\textbf{M})$ implies $\mathcal{A}_{\textbf{M}}(S_{\bgamma})$ is (s) quasi-analytic without the additional assumption~$(\ref{e82})$.
\end{obse}

We now study the quasi-analyticity of the class $\mathcal{A}_{\textbf{M}}(S_{\bgamma})$.

\begin{prop}\label{propcaraccasianalit}
Let $n\in\N$, $\textbf{M}=(M_{p})_{p\in\N_{0}}$ be a sequence, and $\bgamma=(\ga_{1},\ldots,\ga_{n})\in(0,\infty)^{n}$. The following statements are equivalent:
\begin{lst}
\item[(i)] The class $\mathcal{A}_{\textbf{M}}(S_{\bgamma})$ is quasi-analytic.
\item[(ii)] For every $j\in\{1,..,n\}$ the class $\mathcal{A}_{\textbf{M}}(S_{\ga_{j}})$ is quasi-analytic.
\item[(iii)] $\displaystyle\int^{\infty}\frac{\log T_{\widetilde{\textbf{M}}}(r)}{r^{1+1/(\underline{\gamma}+1)}}\,dr$ does not converge.
\end{lst}
If $\textbf{M}$ fulfills $(\a_0)$, the following is also equivalent to the previous statements:
\begin{lst}
\item[(iv)] $\displaystyle\sum_{p=0}^{\infty}\Big(\frac{M_{p}}{(p+1)M_{p+1}}\Big)^{1/(\underline{\gamma}+1)}$ does not converge.
\end{lst}
\end{prop}

\begin{dem}
The equivalence between (iii) and (iv), when $(\a_0)$ is fulfilled, has already been discussed. In order to guarantee the equivalence between (ii) and (iii) we only have to make use of Theorem~\ref{teo30}, and take into account that $\underline{\gamma}\le\ga_{j}$ for every $j\in\mathcal{N}$, that there exists $j$ for which the equality holds, and that, when we replace the value of $\underline{\gamma}$ by a greater one, the integrand increases.

We now prove that (i) implies (ii). Let us suppose there exists $j\in\mathcal{N}$ such that $\mathcal{A}_{\textbf{M}}(S_{\ga_{j}})$ is not quasi-analytic. Consider a non zero map $f_{j}\in\mathcal{A}_{\textbf{M}}(S_{\ga_{j}})$ such that $\mathcal{B}(f_{j})$ is null. The map $f$ defined in $S_{\bgamma}$ by $f(\bz)=f_{j}(z_{j})$, for $\bz=(z_{1},\ldots,z_{n})$, is clearly an element of $\mathcal{A}_{\textbf{M}}(S_{\bgamma})$, it is not null, and $\mathcal{B}(f)=(0)_{\balpha\in\N_{0}^{n}}$: indeed, if $\balpha$ is such that $\balpha_{j'}=\textbf{0}_{j'}$, then
$$\lim_{\bz\to\textbf{0},\bz\in S_{\bgamma}}D^{\balpha}f(\bz)=\lim_{z_j\to 0,z\in S_{\ga_j}}f_{j}^{(\a_{j})}(z_{j})=0,$$
while if $\balpha_{j'}\neq\textbf{0}_{j'}$, then
$$\lim_{\bz\to\textbf{0},\bz\in S_{\bgamma}}D^{\balpha}f(\bz)=\lim_{\bz\to\textbf{0},\bz\in S_{\bgamma}}0=0.$$
We then deduce that $\mathcal{A}_{\textbf{M}}(S_{\bgamma})$ is not  quasi-analytic.\par

Finally, we will prove (ii)$\Rightarrow$(i). Let $n\in\N$, $n>1$, and let $f\in\mathcal{A}_{\textbf{M}}(S_{\bgamma})$ be such that $\mathcal{B}(f)=(0)_{\balpha\in\N_{0}^{n}}$. For every $\balpha_{1'}\in\N_{0}^{1'}$, the map $f_{\balpha_{1'}}\in\textrm{TA}(f)$ belongs to $\mathcal{A}_{\textbf{M}}(S_{\ga_{1}})$ and also $\mathcal{B}(f_{\balpha_{1'}})=(0)_{m\in\N_{0}}$, due to the coherence conditions given in~(\ref{limcondcohe}), since for every $m\in\N_{0}$ we have
$$\lim_{z_{1}\to0,z_{1}\in S_{\ga_1}}f_{\balpha_{1'}}^{(m)}(z_{1})=\lim_{\bz\to\textbf{0},z\in S_{\bgamma}} D^{m\be_1}f(\bz)=0.$$
By virtue of (ii) we have $\mathcal{A}_{\textbf{M}}(S_{\ga_{1}})$ is quasi-analytic, so $f_{\balpha_{1'}}$ is null in $S_{\ga_{1}}$ for every $\balpha_{1'}\in\N_{0}^{1'}$.
Let us fix $z_{1}\in S_{\ga_{1}}$. For every $\balpha_{\{1,2\}'}\in\N_{0}^{\{1,2\}'}$ we consider the map  $f_{\balpha_{\{1,2\}'}}(z_{1},\cdot):S_{\ga_{2}}\to\C$. As we have $f_{\balpha_{\{1,2\}'}}\in\mathcal{A}_{\textbf{M}}(S_{\bgamma_{\{1,2\}}})$
, it is clear that $f_{\balpha_{\{1,2\}'}}(z_{1},\cdot)\in\mathcal{A}_{\textbf{M}}(S_{\ga_{2}})$. Moreover, due to the coherence conditions,
$$
\lim_{z_{2}\to0,z_{2}\in S_{\ga_{2}}} D^{(0,m)}f_{\balpha_{\{1,2\}'}}(z_{1},z_{2})=f_{(m_{\{2\}'},\balpha_{\{1,2\}'})}(z_{1})=0,\quad\hbox{for every }m\in\N_{0}.$$
Therefore, $\mathcal{B}(f_{\balpha_{\{1,2\}'}}(z_{1},\cdot))=(0)_{m\in\N_{0}}$, and as $\mathcal{A}_{\textbf{M}}(S_{\ga_{2}})$ is quasi-analytic, we have $f_{\balpha_{\{1,2\}'}}(z_{1},\cdot)$ is null in $S_{\ga_{2}}$. Varying $z_{1}$ in $S_{\ga_{1}}$ we conclude $f_{\balpha_{\{1,2\}'}}\equiv 0$ in $S_{\ga_{\{1,2\}'}}$ for every $\balpha_{\{1,2\}'}\in\N_{0}^{\{1,2\}'}$.\par
In the next step, only necessary if $n>2$, we fix $\bz_{\{1,2\}}\in S_{\bgamma_{\{1,2\}}}$ and for every $\balpha_{\{1,2,3\}'}\in\N_{0}^{\{1,2,3\}'}$ we can prove, in a similar way, that the map $f_{\balpha_{\{1,2,3\}'}}(\bz_{\{1,2\}},\cdot):S_{\ga_{3}}\to\C$, which belongs to $\mathcal{A}_{\textbf{M}}(S_{\ga_{3}})$, is null. So, we have $f_{\balpha_{\{1,2,3\}'}}\equiv0$ in $S_{\bgamma_{\{1,2,3\}}}$ for every $\balpha_{\{1,2,3\}'}\in\N_{0}^{\{1,2,3\}'}$. After $n$ steps, we conclude that $\textrm{TA}(f)$ is null. Now, condition (ii) together with Proposition~\ref{prop40} guarantees that $\mathcal{A}_{\textbf{M}}(S_{\bgamma})$ is (s) quasi-analytic, so $f$ will identically vanish, as desired.
\end{dem}

We now give a consequence of Proposition~\ref{prop17}.
\begin{prop}\label{propimplicasianalit}
Let $n\in\N$, $\textbf{M}=(M_{p})_{p\in\N_{0}}$ be a sequence of positive real numbers, and $\bgamma=(\ga_{1},\ldots,\ga_{n})\in(0,\infty)^{n}$. If $\displaystyle\int^{\infty}\frac{\log T_{\textbf{M}}(r)}{r^{1+1/\underline{\gamma}}}\,dr$ does not converge, $\mathcal{A}_{\textbf{M}}(S_{\bgamma})$ is  quasi-analytic.
\end{prop}

\begin{dem}
Since $\underline{\ga}\le\ga_j$ for every $j\in\mathcal{N}$, the integrals
$$
\displaystyle\int^{\infty}\frac{\log T_{\textbf{M}}(r)}{r^{1+1/{\gamma_j}}}\,dr
$$
do not converge, so the classes $\mathcal{A}_{\textbf{M}}(S_{\ga_{j}})$, $j\in\mathcal{N}$, are quasi-analytic, and we conclude the proof by the previous result.
\end{dem}

If $\textbf{M}$ is a strongly regular sequence, we can get a consequence of Proposition~\ref{propgammapequecasianal}.

\begin{prop}\label{propgammasubpequecasianal}
Let $n\in\N$, $\textbf{M}=(M_{p})_{p\in\N_{0}}$ be a strongly regular sequence, and $\bgamma=(\ga_1,\ldots,\ga_n)\in(0,\infty)^n$ such that $0<\underline{\ga}<\ga(\textbf{M})$.
Then, the class $\mathcal{A}_{\textbf{M}}(S_{\bgamma})$ is non-quasi-analytic.
\end{prop}

\begin{dem}
There exists $j\in\mathcal{N}$ such that $\ga_j=\underline{\ga}<\ga(\textbf{M})$. According to Proposition~\ref{propgammapequecasianal}, $\mathcal{A}_{\textbf{M}}(S_{\gamma_j})$ is non quasi-analytic. By Proposition~\ref{propcaraccasianalit}, $\mathcal{A}_{\textbf{M}}(S_{\bgamma})$ is also non quasi-analytic.
\end{dem}

Under condition~(\ref{e82}), we can prove the following equivalence, obtaining a new Watson's Lemma type result.

\begin{prop}[Second generalization of Watson's Lemma]
\label{proplemaWatsoncasianalit}
Let $\textbf{M}$ be a strongly regular sequence that verifies condition~(\ref{e82}), and let $n\in\N$ and $\bgamma\in(0,\infty)^{n}$.
The following statements are equivalent:
\begin{lst}
\item[(i)] $\underline{\gamma}\ge\gamma(\textbf{M})$.
\item[(ii)] The class $\mathcal{A}_{\textbf{M}}(S_{\bgamma})$ is quasi-analytic.
\end{lst}
\end{prop}

\begin{dem}
(i) implies (ii) is the only thing left to prove. If (i) is fulfilled, for every $j\in\mathcal{N}$ we have $\gamma_j\ge\gamma(\textbf{M})$. Proposition~\ref{prop18} guarantees that $\mathcal{A}_{\textbf{M}}(S_{\gamma_j})$ is quasi-analytic, and we conclude by applying Proposition~\ref{propcaraccasianalit}.
\end{dem}

\begin{obse}
As it was indicated in Remark~\ref{probabie}, we do not know whether the condition $\underline{\gamma}\ge\gamma(\textbf{M})$ implies $\mathcal{A}_{\textbf{M}}(S_{\bgamma})$ is  quasi-analytic without the additional assumption~$(\ref{e82})$.
\end{obse}

%
%
%
%
%
%
%
%
%

\noindent Dpto. de An\'alisis Matem\'atico y Did\'actica de la Matem\'atica\par\noindent
Facultad de Ciencias, Universidad de Valladolid\par\noindent
Paseo del Prado de la Magdalena s/n\par\noindent
47005 Valladolid, Spain

%
%
%

\end{document}